\documentclass[12pt]{amsart}

\usepackage{amsmath}
\usepackage[utf8]{inputenc}
\usepackage{amssymb}
\usepackage{amsopn}
\usepackage{epsfig}
\usepackage{amsfonts}
\usepackage{latexsym}
\usepackage{enumitem}
\usepackage{color}
\usepackage{graphicx}

\newtheorem{theorem}{Theorem}[section]
\newtheorem{lemma}[theorem]{Lemma}
\newtheorem{proposition}[theorem]{Proposition}

\theoremstyle{definition}

\newtheorem{example}[theorem]{Example}

\theoremstyle{remark}

\numberwithin{equation}{section}

\newcommand{\D}{\ensuremath{\mathcal{D}}}
\newcommand{\N}{\ensuremath{\mathbb{N}}}

\renewcommand{\c}{ {\mathbf{c}}}

\newcommand{\U}{\mathcal{U}}

\newcommand{\set}[1]{\left\{#1\right\}}

\newcommand{\ga}{\gamma}

\newcommand{\f}{\infty}

\newcommand{\om}{\omega}
\newcommand{\al}{\alpha}

\newcommand{\ra}{\rightarrow}

\begin{document}

\title{On small univoque bases of real numbers}

\author{Derong Kong$^1$}
\footnote{The   author  was supported by   NSFC No.~11401516 and Jiangsu Province Natural
Science Foundation for the Youth No.~BK20130433.}
\address{School of Mathematical Science, Yangzhou University, Yangzhou, JiangSu 225002, People's Republic of China.}
\email{derongkong@126.com}

\dedicatory{}


\begin{abstract}
Given a positive real number $x$,  we consider the smallest base $q_s(x)\in(1,2)$ for which there exists a unique sequence $(d_i)$ of zeros and ones such that
\[
x=\sum_{i=1}^\f\frac{d_i}{(q_s(x))^i}.
\]
In this paper we give complete characterizations of those $x$'s for which $q_s(x)\le q_{KL}$,  where $q_{KL}$ is the Komornik-Loreti constant.  Furthermore, we show that $q_s(x)=q_{KL}$ if and only if
\[
x\in\set{1, ~\frac{q_{KL}}{q_{KL}^2-1},~ \frac{1}{q_{KL}^2-1}, ~\frac{1}{q_{KL}(q_{KL}^2-1)}}.
\]
 Finally, we  determine the explicit value of $q_s(x)$ if $q_s(x)<q_{KL}$.
\end{abstract}
\maketitle

\section{Introduction}\label{sec: Introduction}
Given $q\in(1,2)$, a sequence $(d_i)=d_1d_2\cdots$ of zeros and ones is called a \emph{$q$-expansion} of $x$, if
\[
x=\sum_{i=1}^\f\frac{d_i}{q^i}=: ((d_i))_q.
\]
Clearly, a real number $x$ has a $q$-expansion if and only if $x\in I_q:=[0, 1/(q-1)]$.

Non-integer base expansions   pioneered by R\'{e}nyi \cite{Renyi_1957} and Parry \cite{Parry_1960}   obtained great attention from different branches of mathematics, such as number theory, dynamical system, measure theory, combinatorics, et al. In 1990s Erd\H{o}s and Jo\'{o} \cite{Erdos_Joo_1992} discovered that there exist infinitely many reals having a continuum of expansions, and later Sidorov \cite{Sidorov_2003} showed that this property is generic which turns out to be quite different from integer base expansions. Surprisingly, Erd\H{o}s et al. \cite{Erdos_Horvath_Joo_1991, Erdos_Joo_Komornik_1990} also showed that there exist infinitely many reals having a unique expansion. After that there are many works devoted to the investigations of unique expansions (cf.~\cite{Daroczy_Katai_1993, Glendinning_Sidorov_2001, DeVries_Komornik_2008, Komornik_2011, Kong_Li_2015, Komornik_Kong_Li_2015_1}).

On the other hand, let $\U$ be the set of \emph{univoque bases} $q\in(1,2)$ such that $1$ has a unique $q$-expansion. Erd\H{o}s et al. \cite{Erdos_Horvath_Joo_1991} showed that $\U$ is a Lebesgue null set and of first Category. Later, Dar\'{o}czy  and K\'{a}tai \cite{Darczy_Katai_1995} proved that $\U$ has full Hausdorff dimension. Recently, Komornik and Loreti \cite{Komornik_Loreti_2007} investigated the topological properties of $\U$ and showed that its closure $\overline{\U}$ is a Cantor set.

In general, for a real number $x\ge 0$ we consider the set  $\U(x)$   of univoque bases $q\in(1,2)$ such that $x$ has a unique $q$-expansion, i.e.,
\[
\U(x):=\set{q\in(1,2): x~\textrm{has a unique}~q\textrm{-expansion}}.
\]
Clearly, for $x=0$ we have $\U(0)=(1,2)$ since $0$ always has a unique $q$-expansion $0^\f$ for each $q\in(1,2)$. Avoiding this trivial case we will   assume $x>0$ throughout the paper. L\"{u} et al.~\cite{Lu_Tan_Wu_2014} showed that for any $x\in(0,1)$ the set $\U(x)$ is a Lebesgue null set but has full Hausdorff dimension.

When $x=1$, Komornik and Loreti \cite{Komornik_Loreti_1998} considered the smallest base  of $\U(1)$, denoted by $q_{KL}$, which is  called the \emph{Komornik-Loreti constant} in \cite{Glendinning_Sidorov_2001}.  Later, Allouche and Cosnard \cite{Allouche_Cosnard_2000} showed that  $q_{KL}$ is a transcendental number.

In this paper we consider the infinimum  base $q_s(x)$ of $\U(x)$, i.e.,
\[
q_s(x):=\inf\U(x).
\]
Then by \cite{Komornik_Loreti_1998} we have $q_s(1)=q_{KL}\approx 1.78723\in\U(1)$.

Recall that $q_G=(1+\sqrt{5}\,)/2\approx 1.618$ is the golden ratio.
Now we state our main results for $q_s(x)$.

\begin{figure}[h!]
  \centering
  \includegraphics[width=10cm]{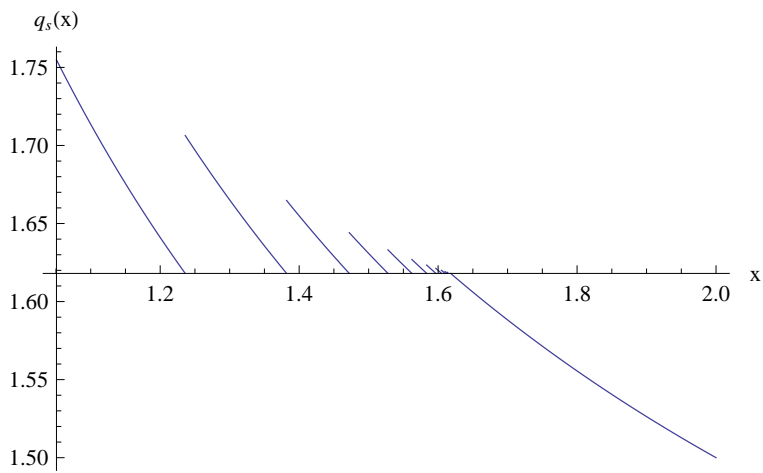}\\
  \caption{Graph of the function $q_s(x)$ with $x\in[z_2, 2]\approx[1.0507, 2]$.}\label{fig:1}
\end{figure}

\begin{theorem}\label{t11}
  Let $x>0$.
Then $q_s(x)\le q_{KL}$ if and only if
\[
\begin{split}
x &\notin  \bigcup_{k=1}^3[(0^{k}(10)^\f)_{q_G}, (0^{k-1}(10)^\f)_{q_{KL}})\\
&\approx [0.236068, 0.255002)\cup[0.381966, 0.455748)\cup[0.618034,0.814527).
\end{split}\]
Furthermore, $q_s(x)=q_{KL}$ if and only if
\[
x\in\set{(0(01)^\f)_{q_{KL}}, ((01)^\f)_{q_{KL}}, ((10)^\f)_{q_{KL}}, 1}.
\]
\end{theorem}
In the following theorem we show that $q_s(x)$ is indeed the smallest base of $\U(x)$ when $q_s(x)\le q_{KL}$.
\begin{theorem}\label{t12}
If $q_s(x)\le q_{KL}$, then  $q_s(x)\in\U(x)$.
  \end{theorem}
We point out that in Theorem \ref{t52} we  determine  the explicit value of $q_s(x)$ when $q_s(x)<q_{KL}$ (see Figure \ref{fig:1} for the graph of $q_s(x)$ with $x\in[1.0507, 2]$).

The rest of the paper is arranged as follows. In Section \ref{sec:unique expansions} we recall some properties of unique expansions.
The proof of Theorem \ref{t11} will be given respectively  in Section \ref{sec: proof of thm11} for the case $x>1$   and  in Section \ref{sec: proof of thm12} for the case $x\in(0,1)$.
In Section \ref{sec: algorithm} we   determine the explicit value of $q_s(x)$ when  $q_s(x)<q_{KL}$,  and prove Theorem \ref{t12}. Finally, we end the paper with some questions.

\section{Unique expansions}\label{sec:unique expansions}
In this section we recall some results of unique expansions.
For $q\in(1,2)$ let $\U_q$ be the set of $x\in I_q$ having a unique $q$-expansion, and let $\U_q'$ be the set of corresponding $q$-expansions.

Recall from \cite{Allouche_Shallit_1999} that $(\tau_i)_{i=0}^\f$ is the classical Thue-Morse sequence beginning with
\[
0\;1\;10\; 1001\; 10010110\; 1001011001101001\;\cdots.
\]
\begin{proposition}\label{p21}
  The classical Thue-Morse sequence $(\tau_i)_{i=0}^\f$ satisfies
  \[
  \tau_0=0,\quad\tau_{2^n}=1\quad\textrm{and}\quad \tau_{2^n+j}=1-\tau_j\quad\textrm{for any}\quad 1\le j<2^n.
  \]
\end{proposition}
For $n\ge 0$ let $q_n\in(1,q_{KL})$ be the appropriate root of the equation
\begin{equation}\label{e21}
1=\sum_{i=1}^{2^n}\frac{\tau_i}{q^i}.
\end{equation}
Then $q_0=1$, $q_1=q_G\approx 1.61803$, $q_2\approx 1.75488$, et al., and $q_n$ strictly increases to $q_{KL}$ as $n\ra\f$.

For $q\in(1,2)$ we denote by $\al(q)=(\al_i(q))$ the quasi-greedy $q$-expansion of $1$, i.e., the lexicographically largest $q$-expansion of $1$ with infinitely many non-zero elements (cf.~\cite{Daroczy_Katai_1993, DeVries_Komornik_2008}).
Then by (\ref{e21}) and Proposition \ref{p21} one can verify that for any $n\ge 1$ we have
\begin{equation}\label{e22}
\al(q_n) =(\tau_1\cdots\tau_{2^{n}}^-)^\f=(\tau_1\cdots\tau_{2^{n-1}}\,\overline{\tau_1\cdots\tau_{2^{n-1}}})^\f.
\end{equation}

Here, for a word $\om=\om_1\cdots\om_m$ we denote by $\om^-:=\om_1\cdots\om_{m-1}(\om_m-1)$ if $\om_m=1$, and we denote by $\om^+:=\om_1\cdots\om_{m-1}(\om_m+1)$ if $\om_m=0$. Furthermore, $\overline{\om}=(1-\om_1)\cdots(1-\om_m)$ stands for the \emph{reflection} of $\om$. For a positive integer $n\ge 1$ we denote by $\om^n$ the concatenations of $\om$ to itself $n$ times,  and   by $\om^\f=\om\om\cdots$ the    concatenations  of $\om$ to itself
infinitely many times.

The following description of the set $\U_q'$ with $q\in(1, q_{KL})$   was essentially established by Glendinning and Sidorov \cite{Glendinning_Sidorov_2001}.
\begin{proposition}
  \label{p22}
  Let $q\in(q_{n-1}, q_{n}]$ with $n\in\N:=\set{ 1, 2, \cdots }$. Then $\U_q'$ contains all sequences of the form
  \[
 (\tau_1^-)^\f,\quad  (\tau_1^-)^*(\tau_1\tau_2^-)^\f,\quad \cdots, \quad(\tau_1^-)^*(\tau_1\tau_2^-)^*\cdots(\tau_1\cdots\tau_{2^{n-1}}^-)^\f,
  \]
  and their reflections, where $*$ stands for all possible non-negative integers.
\end{proposition}
By using Proposition \ref{p22} we can write down the sets $\U_q'$ for $q\in(q_{n-1}, q_{n}]$ and $ n\in\set{1,2, 3}$.
\begin{itemize}
\item If $q\in(q_0, q_1]=(1,(1+\sqrt{5})/2]$, then by Proposition \ref{p22} we have
\[
\U_q'=\set{0^\f, 1^\f}.
\]
\item If $q\in(q_1, q_2]$, then by Proposition \ref{p22} it gives that
\[
\U_q'=\set{0^\f, 0^*(10)^\f, ~1^\f, 1^*(01)^\f}.
\]

\item If $q\in(q_2, q_3]$, then by Proposition \ref{p22} we obtain
\[
\U_q'=\set{0^\f, 0^*(10)^\f, 0^*(10)^*(1100)^\f,~ 1^\f, 1^*(01)^\f, 1^*(01)^*(0011)^\f}.
\]
\end{itemize}

The following monotonicity property of $\U_q'$ is well-known (see, e.g., \cite{Glendinning_Sidorov_2001,Komornik_Loreti_2002, DeVries_Komornik_2008}).
\begin{proposition}
  \label{p23}
  \begin{itemize}
 \item[\rm{(a)}.] $\U_p'\subseteq\U_q'$ for any $1<p<q<2$.
 \item[\rm{(b)}.] Let $q\in(1,2)$ and $(c_i), (d_i)\in\U_q'$. Then $(c_i)<(d_i)$ if and only if $((c_i))_q<((d_i))_q$.
 \end{itemize}
\end{proposition}
 Here and throughout the paper we will use lexicographcial order $<$ or $\le$ between sequences.

\section{Estimation of $q_s(x)$ for  $x>1$}\label{sec: proof of thm11}
In this section we will consider $q_s(x)$ for $x>1$, and prove Theorem \ref{t11} for $x>1$.

For $n\in\N$ we define $\D_n:=\bigcup_{q_{n-1}<q\le q_n}\U_q$. Note by Proposition \ref{p22}   that $\U_q'$ does not change for any $q\in(q_{n-1}, q_n]$. Then
\begin{equation}\label{e31}
\D_n=\bigcup_{q_{n-1}<q\le q_n}\set{((d_i))_q: (d_i)\in\U_{q_n}'}.
\end{equation}

\begin{lemma}
  \label{l31}
Let $N\in\N$. Then  $x\in\bigcup_{n=1}^N\D_n$ if and only if $q_s(x)\le q_N$.
\end{lemma}
\begin{proof}
  If $x\in\bigcup_{n=1}^N \D_n$, then $x\in\D_n$ for some $1\le n\le N$. So, there exist  $q\in(q_{n-1}, q_n]$ and $(d_i)\in\U_{q_n}'$ such that
  \[
  x=((d_i))_q\in\U_q.
  \]
  This implies that $q_s(x)\le q\le q_n\le q_N$.

  On the other hand, suppose that $q_s(x)>q_N$. Then $q>q_N$ for any $q\in\U(x)$. This implies that $x\notin\bigcup_{n=1}^N\D_n$.
\end{proof}
Note that $q_n$ increases to $q_{KL}$ as $n\ra\f$. Then by Lemma \ref{l31} it follows that $q_s(x)<q_{KL}$ if and only if $x\in\bigcup_{n=1}^\f\D_n$. So,
 in order to prove Theorem \ref{t11}  it suffices to show that the union of all $\D_n$  covers the interval $(1,\f)$.

First we construct a sequence $(z_n)$ in $(1,\f)$, where
\[
z_n:=\big(\tau_1\cdots\tau_{2^{n-1}}(\tau_{2^{n-1}+1}\cdots\tau_{2^n})^\f\big)_{q_n}=\left(\tau_1\cdots \tau_{2^{n-1}}(\overline{\tau_1\cdots\tau_{2^{n-1}}}\,^+)^\f\right)_{q_n}.
\]
Here the second equality holds by Proposition \ref{p21}. We will show in Lemma \ref{l33} that $z_n\in\U_{q_n}$ and the sequence $(z_n)$ strictly decreases to $1$.
\begin{lemma}\label{l32}
  For $n\ge 2$ the set $\U_{q_n}'$ contains all sequences of the form
  \[
  \tau_1\cdots\tau_{2^{n-2}}(\overline{\tau_1\cdots\tau_{2^{n-2}}}\,^+)^k(\overline{\tau_1\cdots\tau_{2^{n-1}}}\,^+)^\f,\quad k\in\set{0}\cup\N.
  \]
\end{lemma}
\begin{proof}
  Note by Proposition \ref{p22} that $\U_{q_n}'$ contains  the sequences
  \begin{equation}\label{e32}
  \overline{(\tau_1^-)^2(\tau_1\tau_2^-)\cdots(\tau_1\cdots \tau_{2^{n-3}}^-)(\tau_1\cdots\tau_{2^{n-2}}^-)^k(\tau_1\cdots\tau_{2^{n-1}}^-)^\f},
  \end{equation}
  where $k\in\set{0}\cup\N$. Observe that $\overline{(\tau_1^-)^2}=11=\tau_1\tau_2$, and by Proposition \ref{p21} we have
  \[
  \tau_1\cdots\tau_{2^i}\overline{\tau_1\cdots\tau_{2^i}}\,^+=\tau_1\cdots\tau_{2^{i+1}}\quad\textrm{for any}\quad i\in\N.
  \]
  Therefore, by (\ref{e32}) it follows that $\U_{q_n}'$ contains the sequences
  \[
  \begin{split}
   &\overline{(\tau_1^-)^2(\tau_1\tau_2^-)\cdots(\tau_1\cdots \tau_{2^{n-3}}^-)(\tau_1\cdots\tau_{2^{n-2}}^-)^k(\tau_1\cdots\tau_{2^{n-1}}^-)^\f}\\
   =&\tau_1\tau_2\, \overline{\tau_1\tau_2}\,^+ \cdots\,\overline{\tau_1\cdots\tau_{2^{n-3}}}\,^+ (\overline{\tau_1\cdots\tau_{2^{n-2}}}\,^+)^k (\overline{\tau_1\cdots\tau_{2^{n-1}}}\,^+)^\f\\
   =&\tau_1\tau_2\cdots \tau_{2^{n-2}}\,(\overline{\tau_1\cdots\tau_{2^{n-2}}}\,^+)^k (\overline{\tau_1\cdots\tau_{2^{n-1}}}\,^+)^\f.
  \end{split}\]
\end{proof}

\begin{lemma}
  \label{l33}
  $z_n\in\U_{q_n}$ for any $n\in\N$. Furthermore, $z_n$ strictly decreases to $1$ as $n\ra\f$.
\end{lemma}
\begin{proof}
By taking $k=1$ in Lemma \ref{l32} and using Proposition \ref{p21} it follows that $z_n\in\U_{q_n}$. Then we only need to prove the monotonicity.

  By (\ref{e21}) it follows that $z_n>1$ for any $n\in\N$. Note that $q_n\ra q_{KL}$ as $n\ra\f$. Then
 \[
 \lim_{n\ra\f}z_n=((\tau_i)_{i=1}^\f)_{q_{KL}}=1.
 \]
 So, it suffices to show that $z_{n+1}<z_n$ for any $n\in\N$. Observe that $q_{n+1}>q_n$. Furthermore, by Lemma \ref{l32} we have
 \[
 \tau_1\cdots\tau_{2^{n-1}}\,\overline{\tau_1\cdots\tau_{2^{n-1}}}\,^+(\overline{\tau_1\cdots\tau_{2^n}}\,^+)^\f\in\U_{q_{n+1}}',
 \]
 and in a similar way as in the proof of Lemma \ref{l32} one can verify that $\tau_1\cdots\tau_{2^{n-1}} (\,\overline{\tau_1\cdots\tau_{2^{n-1}}}\,^+)^\f\in\U_{q_{n+1}}'$.
Then by Propositions \ref{p21} and \ref{p23} (b)  it follows that
 \begin{align*}
   z_{n+1}&=\big(\tau_1\cdots\tau_{2^{n}}(\overline{\tau_1\cdots\tau_{2^n}}\,^+)^\f\big)_{q_{n+1}}\\
   &=\big(\tau_1\cdots\tau_{2^{n-1}}\overline{\tau_1\cdots\tau_{2^{n-1}}}\,^+ (\overline{\tau_1\cdots\tau_{2^n}}\,^+)^\f\big)_{q_{n+1}}\\
   &<(\tau_1\cdots\tau_{2^{n-1}}(\overline{\tau_1\cdots\tau_{2^{n-1}}}\,^+)^\f)_{q_{n+1}}\\
   &<(\tau_1\cdots\tau_{2^{n-1}}(\overline{\tau_1\cdots\tau_{2^{n-1}}}\,^+)^\f)_{q_{n}}\;=z_n.
 \end{align*}
\end{proof}

Now we   prove that the union of all $\D_n$ covers $(1,\f)$.
\begin{lemma}
  \label{l34}
$
  \D_1\cap (1,\f)= [z_1, \f).
 $
\end{lemma}
\begin{proof}
By Proposition \ref{p22} it follows  that $\U_{q_1}'=\set{0^\f, 1^\f}$. Then by (\ref{e31}) we have
\[
\begin{split}
\D_1\cap (1,\f) =\set{(1^\f)_q: q\in(q_0,q_1]}&=\set{(1^\f)_q: q\in(1, q_1]}\\
&= [(1^\f)_{q_1}, \f)=[z_1,\f).
\end{split}
\]
\end{proof}

\begin{lemma}\label{l35}
  For $n\ge 2$ we have
  $
  \D_n\cap(1,\f)\supseteq[z_n, z_{n-1}).
  $
\end{lemma}
\begin{proof}
Fix $n\ge 2$. By
  Lemma \ref{l32}  it follows that $\U_{q_n}'$ contains the sequences
  \[
  \c_k:=\tau_1\cdots\tau_{2^{n-2}}(\overline{\tau_1\cdots\tau_{2^{n-2}}}\,^+)^k(\overline{\tau_1\cdots\tau_{2^{n-1}}}\,^+)^\f,\quad k\in\set{0}\cup\N.
  \]
Note   that
 $\c_1<\c_2<\cdots<\c_k<\c_{k+1}<\cdots.$
  Then by Proposition \ref{p23} (b) and Lemma \ref{l32} it follows that
 \begin{equation}\label{eq:33}
 (\c_1)_q<(\c_2)_q<\cdots<(\c_k)_q<(\c_{k+1})_q<\cdots
 \end{equation}
 for any $q\in(q_{n-1}, q_n]$.
Observe by (\ref{e21}) and Proposition \ref{p21}   that
\begin{align*}
(\c_1)_{q_{n}}&=(\tau_1\cdots\tau_{2^{n-2}}\overline{\tau_1\cdots\tau_{2^{n-2}}}\,^+(\overline{\tau_1\cdots\tau_{2^{n-1}}}\,^+)^\f)_{q_{n}}\\
&=( \tau_1\cdots\tau_{2^{n-1}}\,(\overline{\tau_1\cdots\tau_{2^{n-1}}}\,^+)^\f)_{q_{n}}\\
&>(\tau_1\cdots\tau_{2^{n-1}}\, \overline{\tau_1\cdots\tau_{2^{n-1}}}\,^+ 0^\f)_{q_n}=1
\end{align*}
Then by (\ref{e31}) and (\ref{eq:33})     it follows that
\begin{align*}
  \D_n\cap(1,\f)&\supseteq\bigcup_{k=1}^\f\set{(\c_k)_q: q\in(q_{n-1}, q_n]}=\bigcup_{k=1}^\f\big[(\c_k)_{q_n}, (\c_k)_{q_{n-1}}\big).
\end{align*}
In the following we will show that the union on the right indeed equals $[z_n, z_{n-1})$.

 Observe that
 \begin{align*}
 (\c_1)_{q_n}&=(\tau_1\cdots\tau_{2^{n-2}}\,\overline{\tau_1\cdots\tau_{2^{n-2}}}\,^+(\overline{\tau_1\cdots\tau_{2^{n-1}}}\,^+)^\f)_{q_n}\\
 &=(\tau_1\cdots\tau_{2^{n-1}}(\overline{\tau_1\cdots\tau_{2^{n-1}}}\,^+)^\f)_{q_n}=z_n,
 \end{align*}
 and
 $(\c_k)_{q_{n-1}}$  increases to
$
 (\tau_1\cdots\tau_{2^{n-2}} (\overline{\tau_1\cdots\tau_{2^{n-2}}}\,^+)^\f)_{q_{n-1}}=z_{n-1}.
$
Then by (\ref{eq:33}) it suffices to prove that
$   (\c_k)_{q_{n-1}}>(\c_{k+1})_{q_n}$ {for any} $k\ge 1.$

By Proposition \ref{p21} and  (\ref{e21})--(\ref{e22}) it follows that
 \begin{equation}\label{e34}
 \begin{split}
   (\c_k)_{q_{n-1}}&=(\tau_1\cdots\tau_{2^{n-2}} (\overline{\tau_1\cdots\tau_{2^{n-2}}}\,^+)^k (\overline{\tau_1\cdots\tau_{2^{n-1}}}\,^+)^\f)_{q_{n-1}}\\
   &=(\tau_1\cdots\tau_{2^{n-2}}\overline{\tau_1\cdots\tau_{2^{n-2}}}\,^+ 0^\f)_{q_{n-1}}\\
   &~ +   (0^{2^{n-1}} (\overline{\tau_1\cdots\tau_{2^{n-2}}}\,^+)^{k-1}\, (\overline{\tau_1\cdots\tau_{2^{n-2}}} \,\tau_1\cdots\tau_{2^{n-2}})^\f)_{q_{n-1}}\\
   &=1+(0^{2^{n-1}} (\overline{\tau_1\cdots\tau_{2^{n-2}}}\,^+)^{k} 0^\f)_{q_{n-1}}.
 \end{split}
 \end{equation}
 On the other hand, observe that
 \[
  ((\tau_1\cdots\tau_{2^{n-2}}\overline{\tau_1\cdots\tau_{2^{n-2}}})^\f)_{q_n}<((\tau_1\cdots\tau_{2^{n-2}}\overline{\tau_1\cdots\tau_{2^{n-2}}})^\f)_{q_{n-1}}=1.
  \]
 This implies that
 \begin{align*}
   ((\overline{\tau_1\cdots\tau_{2^{n-1}}}\,^+)^\f)_{q_n}&=((\overline{\tau_1\cdots\tau_{2^{n-2}}} \tau_1\cdots\tau_{2^{n-2}})^\f)_{q_n}
   <(\overline{\tau_1\cdots\tau_{2^{n-2}}}\,^+ 0^\f)_{q_n}.
 \end{align*}
 Whence,
 \begin{align*}
   (\c_{k+1})_{q_n}&=(\tau_1\cdots\tau_{2^{n-2}} (\overline{\tau_1\cdots\tau_{2^{n-2}}}\,^+)^{k+1}(\overline{\tau_1\cdots\tau_{2^{n-1}}}\,^+)^\f)_{q_n}\\
   &<(\tau_1\cdots\tau_{2^{n-2}} (\overline{\tau_1\cdots\tau_{2^{n-2}}}\,^+)^{k+2}0^\f)_{q_n}\\
   &=(\tau_1\cdots\tau_{2^{n-2}}\overline{\tau_1\cdots\tau_{2^{n-2}}}\,^+ 0^{k 2^{n-2}}\overline{\tau_1\cdots\tau_{2^{n-2}}}\,^+ 0^\f)_{q_n}\\
   &\hspace{2cm}+(0^{2^{n-1}}(\overline{\tau_1\cdots\tau_{2^{n-2}}}\,^+)^k 0^\f)_{q_n}\\
   &\le(\tau_1\cdots\tau_{2^{n-1}}\overline{\tau_1\cdots\tau_{2^{n-1}}}\,^+ 0^\f)_{q_n}+(0^{2^{n-1}}(\overline{\tau_1\cdots\tau_{2^{n-2}}}\,^+)^k 0^\f)_{q_n}\\
   &=1+(0^{2^{n-1}}(\overline{\tau_1\cdots\tau_{2^{n-2}}}\,^+)^k 0^\f)_{q_n},
 \end{align*}
 where the second inequality follows by using $n\ge 2$ and $k\ge 1$.

Therefore, by (\ref{e34}) we conclude that $(\c_{k+1})_{q_n}<(c_k)_{q_{n-1}}$.
  \end{proof}

\begin{proof}
  [Proof of Theorem \ref{t11} for $x>1$]
  By Lemmas \ref{l34} and \ref{l35} it follows that
  \begin{equation*}
  \bigcup_{n=1}^\f\D_n\supseteq\bigcup_{n=2}^\f[z_n, z_{n-1})\cup[z_1,\f).
  \end{equation*}
  Note by Lemma \ref{l33} that $z_n$ strictly decreases to $1$ as $n\ra\f$. Then
  \[
  \bigcup_{n=1}^\f\D_n\supseteq(1,\f).
  \]
  Hence, by Lemma \ref{l31} we conclude that $q_s(x)<q_{KL}$ for any $x>1$.
\end{proof}

\section{Estimation of $q_s(x)$ for $x\in(0,1)$}\label{sec: proof of thm12}
In this section we will consider $q_s(x)$ for $x\in(0,1)$, and finish the proof of Theorems \ref{t11}.
Recall from (\ref{e31}) that for $n\in\N$ we have
\[
\D_n=\bigcup_{q_{n-1}<q\le q_n}\U_q=\bigcup_{q_{n-1}<q\le q_n}\set{((d_i))_q: (d_i)\in\U_{q_n}'}.
\]
\begin{lemma}
  \label{l41}
  $\D_1\cap(0, 1)=\emptyset$.
\end{lemma}
\begin{proof}
  By Proposition \ref{p22} it follows that
  $
  \U_{q_1}'=\set{0^\f, 1^\f}.
  $
  Furthermore, for any $q\in(q_0, q_1]$ we have
  \[
 (0^\f)_q=0\quad\textrm{and}\quad (1^\f)_q=\frac{1}{q-1}\ge\frac{1}{q_1-1}>1.
  \]
  Therefore, the lemma follows by the definition of $\D_1$.
\end{proof}

\begin{lemma}\label{l42}
  \[
  \bigcup_{n=2}^\f\D_n\cap (0,1)\supseteq\big(0, (0^3(10)^\f)_{q_1}\big)\cup\bigcup_{k=0}^2\big((0^k(10)^\f)_{q_{KL}}, (0^k(10)^\f)_{q_1}\big).
  \]
\end{lemma}
\begin{proof}
  Note by Proposition \ref{p23} (a) that $\U_{q_2}'\subseteq\U_{q_n}'$ for any $n\ge 3$. Then by Proposition \ref{p22} we have $0^k(10)^\f\in\U_{q_n}'$ for any $n\ge 2$ and $k\ge 0$.
  Observe that $q_n$ strictly increases to $q_{KL}$ as $n\ra\f$. Then
  \begin{align*}
    \bigcup_{n=1}^\f\D_n&\supseteq\bigcup_{k=0}^\f\set{(0^k(10)^\f)_q: q\in(q_1, q_{KL})}\\
    &=\bigcup_{k=0}^\f\big((0^k(10)^\f)_{q_{KL}}, (0^k(10)^\f)_{q_1}\big).
  \end{align*}
  Note by Proposition \ref{p23} (b)  that  the sequence $(0^k(10)^\f)_{q_{KL}}$ strictly decreases to $(0^\f)_{q_{KL}}=0$ as $k\ra\f$, and for $k=0$ we have $((10)^\f)_{q_1}=1$. Then
  \begin{align*}
    \bigcup_{n=1}^\f\D_n\cap (0, 1) &\supseteq\bigcup_{k=0}^\f\big((0^k(10)^\f)_{q_{KL}}, (0^k(10)^\f)_{q_1}\big).
  \end{align*}

  So, to complete the proof it suffices to show that
  \begin{equation}
    \label{e41}
    (0^k(10)^\f)_{q_{KL}}<(0^{k+1}(10)^\f)_{q_1}\quad\textrm{for any}\quad k\ge 3.
  \end{equation}
  We will prove this by induction on $k$.

  First we consider $k=3$. Then by (\ref{e21}) and (\ref{e22}) we obtain
  \[
  (0^3(10)^\f)_{q_{KL}}=\frac{1}{q_{KL}^2(q_{KL}^2-1)}<\frac{1}{q_1^4}=(0^310^\f)_{q_1}=(0^4(10)^\f)_{q_1}.
  \]
  Now we assume that (\ref{e41}) holds for some $k\ge 3$. Then
  \[
  (0^{k+1}(10)^\f)_{q_{KL}}=\frac{(0^k(10)^\f)_{q_{KL}}}{q_{KL}}<\frac{(0^{k+1}(10)^\f)_{q_1}}{q_1}=(0^{k+2}(10)^\f)_{q_1}.
  \]
  By induction this establishes (\ref{e41}).
\end{proof}

By Lemmas \ref{l31},  \ref{l41} and \ref{l42} it follows that $q_s(x)<q_{KL}$ for any
\begin{equation*}
x\in(0,1)\setminus\bigcup_{k=1}^3[(0^{k}(10)^\f)_{q_1}, (0^{k-1}(10)^\f)_{q_{KL}}].
\end{equation*}

In the following lemma we   show that the inclusion in Lemma \ref{l42} is indeed an equality.
\begin{lemma}
  \label{l43}
  For any $n\in\N$ we have
  \[\D_n\cap \bigcup_{k=1}^3 [(0^{k}(10)^\f)_{q_1}, (0^{k-1}(10)^\f)_{q_{KL}}]=\emptyset.\]
\end{lemma}
\begin{proof}
  By Lemma \ref{l41} we may assume $n\ge 2$. Suppose on the contrary that there exists
  \begin{equation}\label{e42}
  x=((d_i))_q\in\D_n\cap [(0^{k}(10)^\f)_{q_1}, (0^{k-1}(10)^\f)_{q_{KL}}]
  \end{equation}
   for some $k\in\set{1,2,3}$, where $(d_i)\in\U_{q_n}'$ and $q\in(q_{n-1}, q_n]$.

   By (\ref{e42}) it follows that
  \[
  ((d_i))_q=x\ge(0^k(10)^\f)_{q_1}>(0^k(10)^\f)_q.
  \]
  Note that $0^k(10)^\f\in\U_{q_n}'$. Then by Proposition \ref{p23} (b) it follows that
  \begin{equation}
    \label{e43}
    (d_i)>0^k(10)^\f.
  \end{equation}
  On the other hand, note that
  \[
  ((d_i))_q=x\le(0^{k-1}(10)^\f)_{q_{KL}}<(0^{k-1}(10)^\f)_q.
  \]
Again, by Proposition \ref{p23} (b) we obtain that
  \begin{equation}
    \label{e44}
    (d_i)<0^{k-1}(10)^\f.
  \end{equation}

  Observe by Proposition \ref{p22} that $0^{k-1}(10)^\f$ is the lexicographically smallest sequence in $\U_{q_n}'$ beginning with $0^{k-1}$. Moreover, $0^k(\tau_1\cdots\tau_{2^{n-1}}^-)^\f$  is the lexicographically largest sequence in $\U_{q_n}'$ starting at $0^k$. Hence, by (\ref{e43}) and (\ref{e44})  it follows that
  \[
 0^k(10)^\f< (d_i)\le 0^k(\tau_1\cdots\tau_{2^{n-1}}^-)^\f.
  \]
 Then by (\ref{e22}) this implies
  \begin{align*}
    x =((d_i))_q&\le(0^k(\tau_1\cdots\tau_{2^{n-1}}^-)^\f)_q\\
    &<(0^k(\tau_1\cdots\tau_{2^{n-1}}^-)^\f)_{q_{n-1}}
     =(0^{k-1}\,10^\f)_{q_{n-1}}\\
     &\le (0^{k-1}10^\f)_{q_1}=(0^k(10)^\f)_{q_1},
  \end{align*}
  leading to a contradiction with   (\ref{e42}).
\end{proof}

\begin{proof}[Proof of Theorem \ref{t11} for $0<x<1$.]
By Lemmas \ref{l31} and \ref{l43} it follows that $q_s(x)\ge q_{KL}$ for any $x\in \bigcup_{k=1}^3[(0^{k}(10)^\f)_{q_1}, (0^{k-1}(10)^\f)_{q_{KL}}]$. This yields that $q_s(x)=q_{KL}\in\U(x)$ if
$
x\in\bigcup_{k=0}^2\set{(0^k(10)^\f)_{q_{KL}}}.
$
So, to complete the proof  it suffices to show that $q_s(x)\ne q_{KL}$ for any
$
x\in\bigcup_{k=1}^3 \big[(0^k(10)^\f)_{q_1}, (0^{k-1}(10)^\f)_{q_{KL}} \big).
$

Suppose on the contrary that there exists
\begin{equation}\label{e45}
x=((d_i))_{q_{KL}}\in[(0^k(10)^\f)_{q_1}, (0^{k-1}(10)^\f)_{q_{KL}})
\end{equation}
for some $k\in\set{1,2,3}$, where $(d_i)\in\U_{q_{KL}}'$. Then by Proposition \ref{p23} (b) it gives that
\begin{equation}\label{e46}
(d_i)<0^{k-1}(10)^\f.
\end{equation}
Note by Proposition \ref{p22} that $0^{k-1}(10)^\f$ is the lexicographically smallest sequence in $\U_{q_{KL}}'$ beginning with $0^{k-1}$. Furthermore, any sequence in $\U_{q_{KL}}'$ starting at $0^k$ can not exceed $0^k\tau_1\tau_2\cdots$. Then by (\ref{e46}) it  follows that
\[
(d_i)< 0^k\tau_1\tau_2\cdots.
\]
Therefore,
\begin{align*}
x =((d_i))_{q_{KL}}< (0^k\tau_1\tau_2\cdots)_{q_{KL}}&=(0^{k-1}10^\f)_{q_{KL}}\\
&<(0^{k-1}10^\f)_{q_1}=(0^k(10)^\f)_{q_1},
\end{align*}
leading to a contradiction with (\ref{e45}).
\end{proof}

\section{Explicit value of $q_s(x)$}\label{sec: algorithm}
In this section we will determine  the explicit value of $q_s(x)$ when $q_s(x)<q_{KL}$, and prove Theorem \ref{t12}.
Recall from (\ref{e21}) that $q_0=1, q_1=\frac{1+\sqrt{5}}{2}, \cdots$, and
$q_k$ strictly increases to $q_{KL}$ as $k\ra\f$.  Therefore,
\begin{equation}\label{e51}
(1,q_{KL})=\bigcup_{n=1}^\f(q_{n-1}, q_n],
\end{equation}
where the unions on the right are pairwise disjoint.

 By Lemma \ref{l31} we can deduce directly  the following characterization of those $x$ for which $q_s(x)\in(q_{n-1}, q_n]$.
\begin{lemma}
  \label{l51}
 Let $x>0$ and $n\in\N$. Then $q_s(x)\in(q_{n-1}, q_n]$ if and only if $x\in\D_n\setminus\bigcup_{k=1}^{n-1}\D_k$.
\end{lemma}

In the following theorem we   determine the exact value of $q_s(x)$ when $x\in\bigcup_{n=1}^\f\D_n$.
\begin{theorem}
  \label{t52}
 Let $x\in\D_n\setminus\bigcup_{k=1}^{n-1}\D_k$ with $n\in\N$. Then $q_s(x)\in(q_{n-1}, q_n]$ is the appropriate root of
  \begin{equation}\label{e52}
  x=\sum_{i=1}^\f\frac{\ga_i}{q^i},
  \end{equation}
  where $(\ga_i)$ is the lexicographically smallest sequence in $\U_{q_n}'$ such that $((\ga_i))_{q_{n-1}}>x$.
\end{theorem}
\begin{proof}
 By Lemma \ref{l51} it follows that $q_s(x)\in(q_{n-1}, q_n]$.
 Let $q_*\in(q_{n-1}, q_n]$ satisfy (\ref{e52}). Note that $(\ga_i)\in\U_{q_n}'$. Then by Proposition \ref{p22} it
 follows that $(\ga_i)\in\U_{q_n}'=\U_{q_*}'$. This implies that $x=((\ga_i))_{q_*}$ has a unique $q_*$-expansion, i.e., $q_*\in\U(x)$.
So, $q_s(x)\le q_*$.   In the following it suffices to prove that $q_s(x)\ge  q_*$.

  Take $q\in\U(x)\cap(q_{n-1}, q_n]$, and assume that $(x_i)$ is the unique $q$-expansion of $x$. Then by Proposition \ref{p22} it follows that  $(x_i)\in\U_q'=\U_{q_n}'$.
  Furthermore, $x=((x_i))_q<((x_i))_{q_{n-1}}$.
 By the definition of $(\ga_i)$ we have
 \[(\ga_i)\le (x_i).\] Hence, by using
$
((x_i))_{q}=x=((\ga_i))_{q_*}
 $
 we conclude that $ q\ge q_*$. Therefore, $q_s(x)\ge q_*$.
\end{proof}
Note by (\ref{e51}) that the union of the intervals $(q_{n-1}, q_n]$ for all $n\in\N$ covers $(1, q_{KL})$. Then by Lemma \ref{l51} and Theorem \ref{t52} we can
determine the exact value of $q_s(x)$ whenever  $q_s(x)\in(1, q_{KL})$.

\begin{proof}
  [Proof of Theorem \ref{t12}] Note by Theorem \ref{t11} that $q_s(x)\in\U(x)$ if $q_s(x)=q_{KL}$. Furthermore, if $q_s(x)<q_{KL}$, i.e., $x\in\bigcup_{n=1}^\f\D_n$, then by Theorem \ref{t52} we also have $q_s(x)\in\U(x)$.
\end{proof}

  As an application of Theorem \ref{t52}  we give the explicit value of $q_s(x)$ when $x\ge z_1=1/(q_1-1)$.

 \begin{proposition}
  \label{p53}
  Let $x\ge 1/(q_1-1)$. Then
 $q_s(x)=\frac{1}{x}+1.$
\end{proposition}
\begin{proof}
Take $x\ge 1/(q_1-1)$. Then $x\in\D_1$. Observe that $\U_{q_1}'=\set{0^\f, 1^\f}$. Then by Theorem \ref{t52} it follows that $q_s(x)\in(1, q_1]$ is the root of
\[
x=(1^\f)_q=\frac{1}{q-1}.
\]
This yields that $q_s(x)=\frac{1}{x}+1$.
\end{proof}

Another application of Theorem  \ref{t52} is to consider $q_s(x)$ for $x\in\D_2\cap (1,\f)$. By Lemma \ref{l35} it follows that
\[
\D_2\cap(1,\f)\supseteq [z_2, z_1)=\big[(1^2(01)^\f)_{q_2}, (1^\f)_{q_1}\big).
\]
For $k\ge 1$ let
\[
z_{1,k}:=(1^k(01)^\f)_{q_1}=\frac{1}{q_1}+\cdots+\frac{1}{q_1^k}+\frac{1}{q_1^k(q_1^2-1)}.
\]
Then  $z_{1,k}$ strictly increases to $z_1$ as $k\ra\f$. Observe that
$
z_{1,1}=((10)^\f)_{q_1}=1<z_2.
$
Then the sequence $(z_{1,k})_{k=1}^\f$ forms a partition of   $[z_2, z_1)$.
\begin{lemma}\label{l54}
\[
[z_2, z_1)=[z_2, z_{1,2})\cup\bigcup_{k=2}^\f[z_{1,k}, z_{1,k+1}).
\]
\end{lemma}
\begin{proof}
 Note by (\ref{e21}) and (\ref{e22}) that
  \begin{align*}
 z_{1,1}=1=(11010^\f)_{q_2}<(11(01)^\f)_{q_2} =z_2
  <(11(01)^\f)_{q_1}=z_{1,2}.
  \end{align*}
  Then the lemma follows by observing that $z_{1,k}$ strictly  increases to $z_1$ as $k\ra\f$.
\end{proof}

\begin{proposition}\label{p55}
Let $k\ge 2$. Then for any $x\in[z_{1,k}, z_{1,k+1})$    the smallest base $q_s(x)\in(q_1,q_2)$   is the appropriate root of
\begin{equation}\label{e53}
x=(1^{k+1}(01)^\f)_q=\frac{1}{q}+\cdots+\frac{1}{q^{k+1}}+\frac{1}{q^{k+1}(q^2-1)}.
\end{equation}

Furthermore, for $x\in[z_2, z_{1,2})$ the smallest base $q_s(x)\in(q_1, q_2]$ is the appropriate root of
(\ref{e53}) with $k=1$.
\end{proposition}
\begin{proof}
Fix $k\in\N$ and take $x\in[z_{1,k}, z_{1,k+1})\cap [z_1, z_2)$. Then $x\in\D_2\setminus\D_1$. By Theorem \ref{t52} it follows that
$q_s(x)\in(q_1, q_2]$ is the appropriate root of
\[
x=\sum_{i=1}^\f\frac{\ga_i}{q^i},
\]
where $(\ga_i)$ is the lexicographically smallest sequence in $\U_{q_2}'$ such that $((\ga_i))_{q_1}>x$. So, it suffices to prove that $(\ga_i)=1^{k+1}(01)^\f$.

Observe that
\[
((\ga_i))_{q_s(x)}=x\ge z_{1,k}=(1^k(01)^\f)_{q_1}>(1^k(01)^\f)_{q_s(x)}.
\]
Then by Proposition \ref{p23} (b) it follows that $(\ga_i)>1^k(01)^\f$. By Proposition \ref{p22} the smallest sequence in $\U_{q_s}'$ which is larger than $1^k(01)^\f$ is $1^{k+1}(01)^\f$. Furthermore,
$(1^{k+1}(01)^\f)_{q_1}=z_{1,k+1}>x$. Therefore, by the definition of $(\ga_i)$ we have
$
(\ga_i)=1^{k+1}(01)^\f $
as required.
\end{proof}
By Lemma \ref{l54} and Proposition \ref{p55} we have a complete description of $q_s(x)$ for $x\in[z_2,z_1)$.

\begin{example}\label{ex56}
  By   Propositions \ref{p53} and \ref{p55} we    plot in Figure \ref{fig:1}  the graph of   $q_s(x)$   for $x\in[z_2,z_1)\cup[z_1,2]=[z_2,2]\approx[1.0507, 2]$.

\end{example}

At the end of this section we present some questions.

\begin{itemize}
 \item By Theorem \ref{t11} it follows  that $q_s(x)>q_{KL}$ for any
  \[
x\in\bigcup_{k=1}^3 \big[(0^k(10)^\f)_{q_1}, (0^{k-1}(10)^\f)_{q_{KL}} \big).
\]
Can we determine the exact value of $q_s(x)$ when $q_s(x)>q_{KL}$?

  \item By Theorem \ref{t12} we know that $q_s(x)\in\U(x)$ if $q_s(x)\le q_{KL}$.
   Is it true that $q_s(x)\in\U(x)$ when $q_s(x)>q_{KL}$?
  \end{itemize}


\end{document}